\documentclass[12pt]{amsart}%
\usepackage{amsmath}%
\usepackage{amssymb}%
\usepackage{amscd}%
\usepackage{color}%
\usepackage{enumitem}%
\usepackage{graphicx}%
\usepackage{mathrsfs}%
\usepackage[all]{xy}%
\usepackage{stmaryrd}
\usepackage{newtxtext,newtxmath}%
\usepackage[all, warning]{onlyamsmath}
\usepackage{amsrefs}
\usepackage{bbm}
\newtheorem{thm}{Theorem}[section]
\newtheorem{crl}[thm]{Corollary}

\newtheorem{lmm}[thm]{Lemma}

\newtheorem{prp}[thm]{Proposition}

\newtheorem{con}[thm]{Conjecture}
\theoremstyle{definition}
\newtheorem{dfn}[thm]{Definition}
\theoremstyle{remark}
\newtheorem*{rem}{Remark}
\def\VolumeNo#1{}
\def\YearNo#1{}           
\def\PagesNo#1{}  
\def\R{\mathbb{R}}
\def\bbE{\mathbb{E}}
\def\bbN{\mathbb{N}}
\def\bbZ{\mathbb{Z}}
\def\sI{\mathscr{I}}

\def\La{\Lambda}
\def\cF{{\mathcal{F}}}
\def\s{\eta}
\def\pr{\operatorname{pr}}
\def\unif{{\operatorname{unif}}}
\def\alt{{\operatorname{alt}}}
\def\Map{{\operatorname{Map}}}
\def\Hom{{\operatorname{Hom}}}
\def\loc{{\operatorname{loc}}}
\def\col{{\operatorname{col}}}

\def\diam#1{\operatorname{diam}(#1)}
\def\pair#1{\langle#1\rangle}

\def\abs#1{|#1|}
\def\dabs#1{|\kern-0.5mm|#1|\kern-0.5mm|}
\def\Dabs#1{\Bigl|\kern-0.3mm\Bigl|\,#1\,\Bigr|\kern-0.3mm\Bigr|}
\def\TS{T^*\!S}
\def\vp{\varphi}
\def\ic{\hat\imath}
\def\Ker{\operatorname{Ker}}
\def\C#1{\operatorname{Consv}^\phi_\nu(#1)}
\def\Im{\operatorname{Im}}
\title[Decomposition of Varadhan Type]{A Decomposition Theorem of Varadhan Type for Closed Co-local Forms}
\author[Bannai]{Kenichi Bannai$^{*\diamond}$}\email{bannai@math.keio.ac.jp}
\author[Sasada]{Makiko Sasada$^{\star\diamond}$}\email{sasada@ms.u-tokyo.ac.jp}
\address{${}^*$Department of Mathematics, Faculty of Science and Technology, Keio University, 3-14-1 Hiyoshi, Kouhoku-ku, Yokohama 223-8522, Japan.}
\address{${}^\star$Department of Mathematics, University of Tokyo, 3-8-1 Komaba, Meguro-ku, Tokyo 606-8502, Japan.}
\address{${}^\diamond$Mathematical Science Team, RIKEN Center for Advanced Intelligence Project (AIP),1-4-1 Nihonbashi, Chuo-ku, Tokyo 103-0027, Japan.}
\thanks{Supported by JST CREST Grant Number JPMJCR1913, KAKENHI 18H05233,
and the UTokyo Global Activity Support Program for Young Researchers.}  
\VolumeNo{x}           
\YearNo{202x}           
\PagesNo{000--000}      
\begin{document}
%

\maketitle

\begin{abstract}      
	In our previous article with Yukio Kametani, we investigated the geometric structure underlying a
	large scale interacting system on infinite graphs, via constructing a suitable cohomology
	theory called uniform cohomology, which reflects the geometric property of the
	microscopic model, using a class of functions called the uniform functions.
	In this article, we introduce the \textit{co-local functions}
	on the geometric structure associated to a large scale interacting system.	
	We may define the notion of uniform functions also for co-local functions.
	However, contrary to the functions appearing in our previous article, 
	the co-local functions reflect the stochastic property of the model, namely the
	probability measure on the configuration space.
	We then prove a decomposition theorem of Varadhan type for closed co-local forms.
	The space of co-local functions and forms contain the space of $L^2$-functions and forms.
	 In the last section, 
	 we state a conjecture concerning the decomposition theorem for the $L^2$-case.
\end{abstract}
\setcounter{section}{-1}

\tableofcontents      

\section{Introduction}\label{sec: intro}

In our article \cite{a} with Yukio Kametani, we investigated the geometric structure underlying a
large scale interacting system on infinite graphs, via constructing a suitable cohomology
theory called uniform cohomology, which reflects the geometric property of the
microscopic model.
In this article, as a first step to consider stochastic data on the model,
we will equip the configuration space of our model with a probability measure.
We will then define a class of functions and forms,
which we call the \textit{co-local} functions and forms, reflecting the property of the probability measure.

In our work, when discussing the model of large scale interacting systems, we will use the term
\emph{locale} to refer to the underlying space of the interacting system modeling the space where the dynamics takes place, the
\emph{set of states} to be the single state space of the model 
representing all of the possible states at a single point, and an \emph{interaction}
to be an operator giving the interaction between states at related points
-- the relation between points expressed via the geometric structure of the locale.  These terms will express
the \emph{roles} played by each object within the model, and the corresponding mathematical object 
may differ depending on the context.  In this article as well as \cite{a}, we take the
locale to be any locally finite simple symmetric directed graph $(X,E)$ which is connected.
Here, $X$ denotes the set of vertices of the graph and $E\subset X\times X$ the set of directed edges.
For this article,
we take the set of states to be a \emph{finite set} 
$S$ with a designated element which we call the \textit{base state},
and we let the interaction to be a map
$\phi\colon S\times S\rightarrow S\times S$
such that
for any pair of states $(s_1,s_2)\in S\times S$ satisfying $\phi(s_1,s_2)\neq(s_1,s_2)$,  
we have  
\[
	\ic\circ\phi\circ\ic\circ\phi(s_1,s_2)=(s_1,s_2),
\]
where $\ic\colon S\times S\rightarrow S\times S$
is the bijection obtained by exchanging the components of $S\times S$.
The interaction in our context expresses the possible change of states on vertices
connected by an edge of the locale.

Associated to the data $((X,E),S,\phi)$, 
we define the \textit{configuration} space $S^X$ to be the set
\[
	S^X\coloneqq\prod_{x\in X}S
\]
with the set of \textit{transitions} $\Phi\coloneqq\{ (\s,\eta^e) \mid  \s\in S^X, e\in E\}$,
where for any $e=(o(e),t(e))\in E\subset X\times X$,
the configuration $\eta^e=(\eta^e_x)\in S^X$ is such that $(\eta^e_{o(e)}, \eta^e_{t(e)})=\phi(\eta_{o(e)}, \eta_{t(e)})$
and $\eta^e_x=\eta_x$ if $x\not\eq o(e),t(e)$.
The pair $(S^X, \Phi)$ form a symmetric graph which we call the \emph{configuration space with 
transition structure}.  This graph represents all possible change of the configuration at a single instant.
This structure is independent of the
\textit{transition rate} -- stochastic data which encodes the expected frequency of the transitions.

A typical example of a locale is given by the Euclidean lattice $(\bbZ^d,\bbE^d)$
for $\bbE^d\coloneqq\{ (x,y)\in\bbZ^d\times\bbZ^d\mid |x-y|=1\}$, where $|x|\coloneqq\sum_{i=1}^d|x_i|$
for $x=(x_1,\ldots,x_d)\in\bbZ^d$,
and the simplest nontrivial example of the set of states is $S=\{0,1\}$ with $0$ taken to be 
the base state.  We may interpret the elements $0$ and $1$ respectively to represent
the non-existence or existence of a particle.   An interaction in this case is given by
$\phi(s_1,s_2)=(s_2,s_1)$, expressing the \emph{exclusion process}.
The configuration space $S^X=\{0,1\}^{\bbZ^d}$ expresses
all of the possible configuration of particles on $X$, and the transition structure $\Phi$
expresses the change of configurations arising from the exchange of particles on an edge
of the locale. 
Typical models, such as the \textit{multi-color exclusion process} and the 
\textit{generalized exclusion process} may
be described using this framework.
 See \cite[\S2.1]{a} for other examples of locales, set of states, and interactions.

In our previous article \cite{a}, we investigated the geometric property of $(S^X,\Phi)$ via
a class of functions which we call
the \textit{uniform functions},
and proved a decomposition theorem for shift-invariant closed uniform forms,
which may be interpreted as a uniform variant of the decomposition theorem 
originally proposed by Varadhan \cite{VY97} for
$L^2$-forms, which plays a crucial role in proving
 the hydrodynamic limit for
non-gradient systems.
The uniform functions are certain functions in $C(S^X_*)$, 
where $S^X_*$ denotes the subset of $S^X$ consisting of configurations at base state
outside a finite number of vertices, and $C(S^X_*)$ denotes the $\R$-linear space of $\R$-valued
functions on $S^X_*$.  The definition of the set $C(S^X_*)$ as well as the definition of uniform
functions are purely algebraic, and is independent of the choice of a probability measure on $S^X$.

In this current article, 
we will equip the configuration space $S^X$ with a probability measure $\mu$.
We will then define a certain variant of the space $C(S^X_*)$ which we call the 
\textit{space of co-local functions}, as a projective limit of local functions
via projections defined via the conditional expectations of the probability measure.
The co-local functions are in fact Martingales indexed by the finite subsets of the 
set of vertices of the locale.
As the name suggests, the space of co-local functions is the $\R$-linear dual of the
space of local functions on $S^X$.  In the case that the probability measure $\mu$
is a product measure of a probability measure on $S$,
we will define the notion of \textit{closed co-local forms},
and prove a Varadhan type decomposition theorem for such forms.
The space of co-local functions and forms contain the space of $L^2$-functions and forms.
In  \S\ref{sec: L2}, we consider $L^2$-functions and forms and formulate
a conjecture giving the Varadhan decomposition in the $L^2$-case.

%
%
%
\section{The Space of Co-local Functions}\label{sec: colocal}
%
%
%

Let the notations be as in the \S\ref{sec: intro}.
In this section, we will consider a probability measure on $S^X$ and the conditional 
expectation.  We let  $\cF$ be the set of all subsets of the state space $S$, which is a
finite set since we have assumed that $S$ is finite.
Then $\cF$ is trivially a $\sigma$-algebra.   In other words, $\cF$ contains $S$ and is closed under
taking complement and countable unions.  
Note that $\cF$ coincides with the Borel $\sigma$-algebra for the discrete topology of $S$.
The condition that $X$ is connected and locally finite ensures that the set of vertices is countable.
For any $\La\subset X$, we let $\cF_\La\coloneqq\cF^{\otimes \La}$ be the 
product $\sigma$-algebra on $S^\La$ obtained from $\cF$,
which coincides with the Borel $\sigma$-algebra for the topological space $S^\La$.
In particular, if $\La$ is finite, then $\cF_\La$ coincides with the set of all subsets of $S^\La$.
For any $\La\subset X$, the pair $(S^\La,\cF_\La)$ is a measurable space.

\begin{dfn} 
	For any $\La\subset X$, denote by $C(S^\La)$ 
	the $\R$-linear space of \textit{real valued measurable functions} on $S^\La$.
\end{dfn}

In particular, if $\La$ is finite, then $C(S^\La)$ is simply the set of real valued functions on $S^\La$.
Any inclusion $\La\subset \La'$ of subsets in $X$,
the inclusion induces a projection $S^{\La'}\rightarrow S^\La$, which in
turn induces the inclusion $C(S^\La)\hookrightarrow C(S^{\La'})$.
We let $\sI$ be the set of \emph{all finite subsets} of $X$.
Then $\sI$ is a \textit{directed set} for the order given by the inclusion.
Indeed, $\subset$ gives a \textit{partial order}, and
for any $\La,\La'\in\sI$, if we let $\La''\coloneqq\La\cup\La'$, then we have $\La,\La'\subset\La''$.
Note that if $\La\subset\La'$, 
then there exists a natural projection $S^{\La'}\rightarrow S^{\La}$,
which induces a natural inclusion $C(S^{\La})\hookrightarrow C(S^{\La'})$.
This gives $\{ C(S^\La)\}_{\La\in\sI}$ a structure of a \textit{directed system}.
\begin{dfn}
	We define the set of \textit{local functions} on $S^X$ to be the direct limit
	\[
		C_\loc(S^X)\coloneqq\varinjlim_{\La\in\sI}C(S^\La)=\bigcup_{\La\in\sI}C(S^\La).
	\]
	This definition of local functions is exactly the definition given in \cite{a}.
\end{dfn}

Let $S^X_*$ be the set of configurations whose components are at base state except for 
a finite number of vertices, and let $C(S^X_*)$ be the $\R$-linear space of real valued functions on $S^X_*$.
In \cite{a}, we defined the space of uniform functions to be a certain subspace of $C(S^X_*)$
containing the space of local functions $C_\loc(S^X)$.  
The space $C(S^X_*)$ may be interpreted as a projective limit of space of 
local functions, given as follows.
For any inclusion $\La\subset\La'\subset X$, consider the projection
$\iota^\La\colon C(S^{\La'})\rightarrow C(S^\La)$ defined by
\[
	\iota^\La f(\s)\coloneqq f(\iota_{\La}(\s))
\]
for any $\s\in S^\La$, where $\iota_{\La}(\s)$ is the configuration $\s'=(\eta'_x)\in S^{\La'}$
such that $\eta'_x=\eta_x$ for any $x\in \La$ and $\eta'_x$ is at base state for any $x\in\La'\setminus\La$.
Then the $\R$-linear spaces $C(S^\La)$ for $\La\in\sI$ form a projective system 
with respect to the projection $\iota^\La$.

\begin{lmm}\label{lem: pro-local}
	We have
	\[
		C(S^X_*)=\varprojlim_{\La\in\sI}C(S^\La).
	\]
\end{lmm}

\begin{proof}
	Suppose $(f^\La)\in \varprojlim_{\La\in\sI}C(S^\La)$.
	For any $\s\in S^X_*$, let $\La\in\sI$ be sufficiently large such that the components $\eta_x$
	of $\s=(\eta_x)$	is at base state outside $x\in\La$.
	Then $f(\s)\coloneqq \iota^\La f(\s)$ is independent of the choice of such $\La$,
	hence this construction defines a function $f\colon S^X_*\rightarrow\R$.
	This gives an $\R$-linear homomorphism 
	\begin{equation}\label{eq: 1}
		\varprojlim_{\La\in\sI}C(S^\La)\rightarrow C(S^X_*).
	\end{equation}
	On the other hand, suppose we have a function $f\colon S^X_*\rightarrow\R$.
	For any $\La\in\sI$, let $f^\La(\s)\coloneqq f(\s|_\La)$ for any $\s\in S^X$,
	where $\s|_\La$ denotes the configuration in $S^X$ whose components at $x\in\La$
	coincides with that of $\s$ and is at base state outside $\La$.
	Then $f^\La$ is a function in $C(S^\La)$, and $(f^\La)$ for $\La\in\sI$ form a projective
	system with respect to the projection $\iota^\La$.
	This gives an $\R$-linear homomorphism inverse to that of \eqref{eq: 1},
	which proves that \eqref{eq: 1} is an isomorphism as desired.
\end{proof}

Due to Lemma \ref{lem: pro-local}, we refer to a function 
$f\in C(S^X_*)$ as a \emph{pro-local function}.
The projection $\iota^\La\colon C(S^X)\rightarrow C(S^\La)$ maps any $f\colon S^X\rightarrow\R$ to a function $\iota^\La f\colon S^\La\rightarrow\R$.
This operation may be interpreted as restricting a function $f$ on $S^X$
to a function which depends only on the local configurations
$S^\La$ on the vertices $\La\subset X$.  One drawback of the projection $\iota^\La$ is 
that by construction, $\iota^\La f$ depends only on the restriction of $f$ to $S^X_*$ and 
does not reflect the behavior of $f$ on the entirety of the configuration space $S^X$.
In order to redress the projection, we introduce a probability measure $\mu$ on $S^X$
encoding the probability of occurrence of the configurations in $S^X$.

We fix a probability measure $\mu$ on $(S^X,\cF_X)$.
For any $\La\subset X$, we denote again by $\mu$ the measure on $S^\La$ 
obtained as the pushforward of the measure $\mu$ with respect to the projection $S^X\rightarrow S^\La$.  We assume in addition that $\mu$ is supported on $S^\La$ for any finite $\La\subset X$.
In other words, we have $\mu(\s)\coloneqq\mu(\{\s\})>0$ for any $\s\in S^\La$.
For any \emph{integrable} 
$f\in C(S^X)$, we let $E_\mu[f]\coloneqq\int_{S^X}fd\mu$ be the expectation value of $f$ on $S^X$
with respect to the measure $\mu$.
For any $\La\subset X$ and an integrable function $f\in C(S^{\La'})$ for $\La'\supset\La$,
we let $\pi^\La_\mu f \coloneqq E_\mu[f|\cF_\La]$ be the \textit{conditional expectation}  
with respect to the projection $\pr_\La\colon S^{\La'}\rightarrow S^\La$.
More precisely, we define $\pi^\La_\mu f\in C(S^\La)$ 
to be the integrable function on $S^\La$ characterized by the property that
\[
	E_\mu[(\pi^\La_\mu f)g]\coloneqq\int_{S^{\La}}(\pi^\La_\mu f) g d\mu = \int_{S^{\La'}}f g d\mu 
\]
for any integrable $g\in C(S^{\La})$,
where we denote again by $g$ the function in $C(S^{\La'})$
induced by the natural inclusion $C(S^\La)\hookrightarrow C(S^{\La'})$.
In particular, if $\La,\La'\in\sI$, 
by taking $g$ to be the \textit{indicator function} $1_\s$
for $\s\in S^\La$
which is \textit{one} on $\s$ and \textit{zero} outside of $\s$,
we see that $\pi^\La_\mu f(\s)$ is explicitly given as
\begin{equation}\label{eq: conditional expectation}
	\pi^\La_\mu f(\s)=\frac{1}{\mu(\s)}\sum_{\substack{\s'\in S^{\La'}\\ \pr_\La(\s')=\s}} f(\s')\mu(\s').
\end{equation}
In other words, $\pi^\La_\mu f(\s)$ is the expected value of the function $f$ on the set of 
configurations in $S^{\La'}$ which projects to $\s\in S^\La$.
Hence if $f$ is an integrable function in $C(S^X)$, then 
$\pi^\La_\mu f$ is a function in $C(S^\La)$ which reflects the property of $f$
on the entirety of $S^X$ weighted by the probability measure $\mu$.
In the case that $f\in C(S^\La)$, if we view $f$ as an element in $C(S^{\La'})$ for $\La\subset\La'$
through
the natural inclusion $C(S^\La)\hookrightarrow C(S^{\La'})$, then 
we see that $\pi_\mu^\La f=f$.

The conditional expectation satisfies the \textit{tower property}, which
may be proved for the case $\La,\La',\La''\in\sI$ such that
$\La\subset \La'\subset \La''$ 
and integrable $f\in C(S^{\La''})$ by
\begin{equation}\label{eq: tower}\begin{split}
	\pi^\La_\mu(\pi^{\La'}_\mu f)(\s)
	&=\frac{1}{\mu(\s)}\sum_{\substack{\s'\in S^{\La'}\\ \pr_\La(\s')=\s}} \pi_\mu^{\La'}f(\s')\mu(\s')
	=\frac{1}{\mu(\s)}\sum_{\substack{\s'\in S^{\La'}\\ \pr_\La(\s')=\s}} 
	\sum_{\substack{\s''\in S^{\La''}\\ \pr_{\La'}(\s'')=\s'}}f(\s'')\mu(s'')\\
	&=\frac{1}{\mu(\s)}\sum_{\substack{\s''\in S^{\La''}\\ \pr_\La(\s'')=\s}} f(\s'')\mu(\s'')
	=\pi^\La_\mu f(\s).
\end{split}\end{equation}
Note that if $\La\in\sI$, then any function in $C(S^\La)$ is integrable.
Hence for any $\La,\La'\in\sI$ such that $\La\subset\La'$, the conditional expectation
induces a homomorphism $\pi^\La_\mu\colon C(S^{\La'})\rightarrow C(S^\La)$,
which is a projection since $\pi^\La_\mu (\pi^\La_\mu f)=\pi^\La_\mu f$
for any $f\in C(S^{\La'})$.
The tower property shows that
$\{C(S^\La)\}_{\La\in\sI}$ form a projective system with respect to the projections
$\pi^\La_\mu$.

\begin{dfn}
	We define the set of \textit{co-local functions} on $S^X$ to be the projective limit
	\[
		C_{\col}(S^X_\mu)\coloneqq\varprojlim_{\La\in\sI}C(S^\La)
	\]
	of $C(S^\La)$ with respect to the projections $\pi^\La_\mu$.
	We call any element $f$ in $C_{\col}(S^X_\mu)$ a \textit{co-local function} on $S^X$.
\end{dfn}

By definition, a co-local function $(f^\La)\in C_{\col}(S^X_\mu)$ is a system of measurable functions 
$f^\La\in C(S^\La)$
for $\La\in\sI$ satisfying $\pi^\La_\mu f^{\La'}=f^\La$ for any $\La,\La'\in\sI$ such that $\La\subset\La'$.
Such a system of random variables related via the conditional expectation is usually referred to as a \textit{Martingale} with respect to the index set $\sI$.
We note that if $X$ is an infinite locale, then 
co-local functions do not necessarily define a function on $S^X$.
We believe the space of co-local functions is a natural framework to
consider formal infinite sum of functions appearing in the works of Varadhan
(see \eqref{eq: CQ} below for the case of conserved quantities).

If $f$ is a function in $C(S^\La)$
for some $\La\in\sI$, then $f$ is an integrable function in $C(S^X)$
through the natural inclusion $C(S^\La)\hookrightarrow C(S^X)$.
By \eqref{eq: conditional expectation}, we see that $\pi^\La_\mu f=f$
as a function in $C(S^\La)$.
For any local function $f\in C_\loc(S^X)$,
if we let $f^\La\coloneqq\pi^\La_\mu f$ for any $\La\in \sI$, then the system
$(f^\La)$ is a co-local function on $S^X$.  Thus we have a homomorphism
\[
	C_\loc(S^X)\rightarrow C_{\col}(S^X_\mu),\qquad f\mapsto (\pi^\La_\mu f)
\]
which is injective since $f=\pi^\La_\mu f$ for $\La$ sufficiently large.
Hence we may view the space of co-local functions as enlarging the space of local functions.

\begin{rem}
	In what follows, we will 
	denote $\pi^\La_\mu$ and $C_{\col}(S^X_\mu)$ simply as $\pi^\La$ and $C_\col(S^X)$.
	If $X$ is a finite locale, then we simply have $C_\loc(S^X)=C_\col(S^X)=C(S^X)$.
	As explained below, $C(S^\La)$ for any $\La\in\sI$ has a structure of a Hilbert space
	with respect to the inner product \eqref{eq: inner product}.
	Since $C_\col(S^X)$ is a projective limit of Hilbert spaces, it has a structure of
	a Fr\'echet space (see for example \cite[Proposition 2.3.7]{DGV16}).
	The space of local functions $C_\loc(S^X)$ is dense in $C_\col(S^X)$ for this topology.
\end{rem}

The term \textit{co-local} 
is derived from the duality between local and co-local functions, given in Lemma \ref{lem: 2}
below.
For any $\La\in\sI$, the space $C(S^\La)$ 
is a Hilbert space with respect to the inner product 
$\pair{\cdot,\cdot}_\mu$ given for any $f,g\in C(S^\La)$ by 
\begin{equation}\label{eq: inner product}
	\pair{f,g}_\mu\coloneqq E_\mu[fg]=\int_{S^\La}fgd\mu.
\end{equation}
The Riesz representation theorem gives an isomorphism
\[
	C(S^\La)\xrightarrow\cong \bigl(C(S^\La)\bigr)^* \qquad f \mapsto 
	(g\mapsto\pair{f,g}_\mu),
\]
where $\bigl(C(S^\La)\bigr)^*$ denotes the 
space of bounded, or equivalently continuous $\R$-linear functionals on $C(S^\La)$.
Since we have assumed that $S$ hence $S^\La$ is finite, 
$C(S^\La)$ is a finite dimensional $\R$-linear space.
Hence any $\R$-linear functional
on $C(S^\La)$ is automatically bounded.
This shows that we have $\bigl(C(S^\La)\bigr)^*= \bigl(C(S^\La)\bigr)^\vee$,
where $\bigl(C(S^\La)\bigr)^\vee$ denotes the algebraic dual 
$\bigl(C(S^\La)\bigr)^\vee\coloneqq\Hom_\R(C(S^\La),\R)$.


\begin{lmm}\label{lem: 2}
	For any $(f^\La)\in C_\col(S^X)$, consider the bounded linear functional on $C(S^\La)$
	given by $g\mapsto \langle f^\La,g \rangle_\mu$.  Then this gives an isomorphism
	\[
		C_\col(S^X) \cong \bigl(C_\loc(S^X)\bigr)^\vee.
	\]
\end{lmm}

\begin{proof}
	For any $\La,\La'\in\sI$ such that $\La\subset\La'$, consider
	$f\in C(S^{\La'})$ and $g\in C(S^\La)$.
	By definition of the conditional expectation,
	we have
	\[
		\pair{\pi^\La f,g}_\mu=\int_{S^\La} (\pi^\La f)g d\mu
		=\int_{S^{\La'}} f  g d\mu
		=\pair{f, g}_\mu.
	\]
		Thus we have a commutative diagram
	\[
	\xymatrix{
		C(S^{\La'})\ar[r]^\cong\ar[d]_{\pi^\La}& C(S^{\La'})^\vee\ar[d]^{\iota^*}\\
		C(S^\La)\ar[r]^\cong& C(S^\La)^\vee,
	}
	\]
	where $\iota^*$ is the dual of the natural injection $\iota\colon C(S^\La)\hookrightarrow C(S^{\La'})$.
	By passing to the projective limit, we have
	\[
		C_\col(S^X)=\varprojlim_{\La\in\sI}C(S^\La)\cong\varprojlim_{\La\in\sI}C(S^\La)^\vee
		=\bigl(\varinjlim_{\La\in\sI}C(S^\La)\bigr)^\vee=(C_\loc(S^X))^\vee
	\]
	as desired.
\end{proof}

%
%
%
\section{Uniformly Co-Local Functions}\label{sec: uniform}
%
%
%

In this section, we give an expansion of co-local functions and define the notion of uniformity.
This is a co-local version of \cite[Proposition 3.3]{a}.
We first prove some properties concerning the conditional expectation.

\begin{lmm}\label{lem: orthogonal}
	Consider $\La,\La'\in\sI$ such that $\La\subset\La'$.
	Then $\pi^\La$ is an orthogonal projection of $C(S^{\La'})$ to $C(S^\La)$ 
	with respect to the inner product \eqref{eq: inner product} on $C(S^{\La'})$.
\end{lmm}

\begin{proof}
	The homomorphism $\pi^\La\colon C(S^{\La'})\rightarrow C(S^\La)$ is a projection
	since $\pi^\La(\pi^\La f)=\pi^\La f$ for any $f\in C(S^{\La'})$.
	For any $f,g\in C(S^{\La'})$, the characterization of the conditional expectation shows
	that
	\begin{align*}
		\pair{\pi^\La f, \pi^\La g}_\mu&=E_\mu[(\pi^\La f)(\pi^\La g)]
		=E_\mu[(\pi^\La f)g]=\pair{\pi^\La f, g}_\mu,\\
		\pair{\pi^\La f, \pi^\La g}_\mu&=E_\mu[(\pi^\La f)(\pi^\La g)]
		=E_\mu[f(\pi^\La g)]=\pair{f,\pi^\La g}_\mu,
	\end{align*}
	which shows that $\pair{\pi^\La f, g}_\mu=\pair{f,\pi^\La g}_\mu$ as desired.
\end{proof}

We now assume until the end of this section
that the probability measure $\mu$ on $S^X$ is the 
product measure $\mu=\nu^{\otimes X}$ for a probability measures $\nu$ on $S$ 
supported on points, i.e.\
satisfying 
$\nu(\{s\})>0$ for any $s\in S$.
In this case, for any $\La\subset X$, the push-forward of $\mu$ to $S^\La$
coincides with the product measure $\nu^{\otimes\La}$.

\begin{lmm}\label{lem: needed}
	Suppose $\mu$ is the product measure $\mu=\nu^{\otimes X}$ as above.
	For any $\La\in\sI$ and $\La',\La''\subset\La$, we have
	\[
		\pi^{\La'}(\pi^{\La''} f) = \pi^{\La'\cap\La''}f
	\]
	for any function $f\in C(S^\La)$.
\end{lmm}

\begin{proof}
	By \eqref{eq: conditional expectation}, for any $\s'\in S^{\La'}$, we have
	\begin{align*}
		\pi^{\La'}f(\s')
		&\coloneqq\frac{1}{\mu(\s')} \sum_{\substack{\s\in S^\La\\\pr_{\La'}(\s)=\s'}} f(\s)\mu(\s)
		=\sum_{\substack{\s=(\eta_x)\in S^\La\\\pr_{\La'}(\s)=\s'}} f(\s)\Bigl(\prod_{x\in\La\setminus\La'}\nu(\eta_x)\Bigr),
	\end{align*}
	where the second equality is derived from the fact that $\mu=\nu^{\otimes X}$.
	Again by \eqref{eq: conditional expectation}, we have
	\begin{align*}
		(\pi^{\La'}(\pi^{\La''}f))(\s')
		&=\sum_{\substack{\s=(\eta_x)\in S^\La\\\pr_{\La'}(\s)=\s'}}\pi^{\La''}f(\s)
		\Bigl(\prod_{x\in\La\setminus\La'}\nu(\eta_x)\Bigr)\\
		&=\sum_{\substack{\s''=(\eta''_x)\in S^{\La''}\\\pr_{\La'\cap\La''}(\s'')=\pr_{\La'\cap\La''}(\s')}}
		\pi^{\La''}f(\s'') \Bigl(\prod_{x\in\La''\setminus\La'}\nu(\eta''_x)\Bigr)\\
		&=\sum_{\substack{\s''\in S^{\La''}\\\pr_{\La'\cap\La''}(\s'')=\pr_{\La'\cap\La''}(\s')}}
		\sum_{\substack{\s\in S^\La\\\pr_{\La''}(\s)=\s''}} f(\s)
		\Bigl(\prod_{x\in\La\setminus\La''}\nu(\eta_x)\Bigr)\Bigl(
		\prod_{x\in\La''\setminus\La'}\nu(\eta''_x)\Bigr),
	\end{align*}
	which coincides with
	\begin{align*}
		\pi^{\La'\cap\La''}f(\s')
		= \sum_{\substack{\s\in S^\La\\\pr_{\La'\cap\La''}(\s)=\s'}} f(\s)
		\Bigl(\prod_{x\in\La\setminus(\La'\cap\La'')}\nu(\eta_x)\Bigr)
	\end{align*}
	as desired.
\end{proof}

\begin{prp}\label{prop: expansion}
	Suppose $\mu$ is the product measure $\mu=\nu^{\otimes X}$ on $S^X$, 
	where $\nu$ is a probability measure
	on $S$ supported on the points of $S$.
	For any $\La\in\sI$, let
	\[
		C_{\La}(S^X)\coloneqq\{f\in C(S^\La)\mid\pi^{\La'} f\equiv0\text{ if }\La\not\subset\La'\}.
	\]
	Then for any $(f^\La)\in C_\col(S^X)$, there exists a unique family of functions $f_\La\in C_{\La}(S^X)$
	such that
	\begin{equation}\label{eq: satisfy}
		f^\La=\sum_{\La''\subset\La} f_{\La''}
	\end{equation}
	for any $\La\in\sI$.
\end{prp}

\begin{proof}
	We prove our result by induction on the order of $\La$.
	We first let $f_{\emptyset}\coloneqq f^\emptyset$, which is the case for $\La=\emptyset$.
	Next, for any $\La\in\sI$, suppose $f_{\La''}$ is defined for any $\La''\subsetneq\La$.
	We let
	\begin{equation}\label{eq: construction}
		f_{\La}\coloneqq f^\La-\sum_{\La''\subsetneq\La} f_{\La''}.
	\end{equation}
	Then for any $\La'\in\sI$ such that $\La\not\subset\La'$, we have
	\[
		\pi^{\La'} f_\La=\pi^{\La'}f^\La - \sum_{\La''\subsetneq\La} \pi^{\La'} f_{\La''}
		=f^{\La\cap\La'}-\sum_{\La''\subset\La\cap\La'} f_{\La''}\equiv 0,
	\]
	where the last equality is from the induction hypothesis.
	In the calculation, we have used the fact that 
	$\pi^{\La'}f^\La=f^{\La\cap\La'}$
	and
	\[
		\pi^{\La'} f_{\La''}
		=\pi^{\La' \cap \La''} f_{\La''}=
		\begin{cases}
			f_{\La''}  &  \La''\subset\La\cap\La' \\
			0 &  \La''\not\subset\La\cap\La',
		\end{cases}
	\]
	which follows from Lemma \ref{lem: needed}.
	This proves that we have $f_{\La}\in C_\La(S^X)$.
	Hence by induction, there exists unique $f_{\La}\in C_\La(S^X)$ for any $\La\in\sI$ satisfying
	\eqref{eq: satisfy} as desired.
\end{proof}

\begin{rem}
	By abuse of notation, we will often write the equality \eqref{eq: satisfy}  as
	\begin{equation}\label{eq: expansion2}
		f=\sum_{\La''\in\sI} f_{\La''}
	\end{equation}
	for $f=(f^\La)_{\La\in\sI}\in C_\col(S^X)$. If $f\in C(S^\La)$ for some $\La\in\sI$,
	then we have $f=\sum_{\La''\subset\La} f_{\La''}$.
\end{rem}

We define the notion of uniformity for co-local functions as follows.
For any $x,x'\in X$, let $d_X(x,x')$ be the length of
the shortest path from $x$ to $x'$ in $X$, and for any
$\La\subset X$, we define the diameter of $\La$ by
$\diam\La\coloneqq\sup_{x,x'\in\La}d_X(x,x')$.

\begin{dfn}\label{def: uniform}
	We say that a co-local function $f$ is \textit{uniformly co-local},
	if there exists $R>0$ such that $f_\La=0$ for any $\La\in\sI$ with $\diam\La>R$
	in the expansion 
	\[
		f=\sum_{\La\in\sI}f_\La
	\]
	 of Proposition \ref{prop: expansion}.
	We denote by $C_\unif(S^X)$ the set of all uniformly co-local functions.	
\end{dfn}
We define the subspaces $C^0_\col(S^X)$ and $C^0_\unif(S^X)$ of $C_\col(S^X)$ and $C_\unif(S^X)$ by
\begin{align*}
		C^0_\col(S^X)&\coloneqq\bigl\{f\in C_\col(S^X) \mid f^\emptyset\equiv0\bigr\},&
		C^0_\unif(S^X)&\coloneqq\bigl\{f\in C_\unif(S^X) \mid f^\emptyset\equiv0\bigr\}.
\end{align*}
Furthermore, for any $\La\in\sI$, we let
\[
	C^0(S^\La)\coloneqq\bigl\{f\in C(S^\La) \mid E_\mu[f]\equiv0\bigr\}.
\]
Since $S^\emptyset$ is a set consisting of a single point, for any local function $f\in C(S^\La)$, 
the function $\pi^\emptyset f$
is the constant functions with value $E_\mu[f]$.  The condition $E_\mu[f]=0$
is equivalent to the condition $\pi^\emptyset f\equiv0$.

The most important example of a uniformly co-local function is given by the 
conserved quantities.

\begin{dfn}\label{def: 10}
	We say that a map $\xi\colon S\rightarrow\R$ is a \textit{$\nu$-regularized conserved quantity} for the
	interaction $\phi$, if $E_\nu[\xi]=0$ and
	\[
		\xi(s_1')+\xi(s_2')=\xi(s_1)+\xi(s_2)
	\]
	for any $(s_1,s_2)\in S\times S$, where $(s_1',s_2')=\phi(s_1,s_2)$.
	We denote by $\C{S}$ the $\R$-linear space of $\nu$-regularized conserved quantities on $S$.
\end{dfn}

The definition of conserved quantities in Definition \ref{def: 10} slightly
differs from that of \cite{a} since we normalize with the condition $E_\nu[\xi]=0$
instead of the condition that $\xi$ is \textit{zero} at the base state.
For any $x\in X$, the natural projection $\pr_{\{x\}}\colon S^X\rightarrow S^{\{x\}}=S$ induces an
inclusion $C(S)\hookrightarrow C(S^X)$.  For any $\xi\in \C{S}$, we denote by $\xi_x$ the image of
$\xi$ with respect to this inclusion, and we let $\xi^\La\coloneqq\sum_{x\in \La}\xi_x$.
The system of local functions $(\xi^\La)$ form a co-local function on $S^X$ which we denote by $\xi_X$.
Then the expansion of $\xi_X$ in Proposition \ref{prop: expansion} is simply
\begin{equation}\label{eq: CQ}
	\xi_X=\sum_{x\in X}\xi_x,
\end{equation}
which shows that $\xi_X$ is uniformly co-local for the constant $R=1$.
Due to our normalization $E_\nu[\xi]=0$, we see that $\xi_X$ is an element in $C^0_\unif(S^X)$.
By associating to $\xi\in\C{S}$ the uniformly co-local function $\xi_X$ in $C^0_\unif(S^X)$,
we have an $\R$-linear homomorphism
\[
	\C{S}\rightarrow C^0_\unif(S^X),
\]
which is injective since $\xi_X$ is \textit{zero} if and only if $\xi_x$ is \textit{zero}
for any $x\in X$, the last condition equivalent to the condition that $\xi$ is constantly \textit{zero}.

%
%
%
\section{Local and Co-local Forms}\label{sec: forms}
%
%
%

In this section, we define the space of local and co-local forms.
We let $\mu$ be a probability measure on $S^X$ which is supported on 
$S^\La$ for any finite $\La\subset X$.
For any $\La\subset X$,
we let 
\[
	\Phi_\La\coloneqq \{  (\s,\s^e)\mid\s\in S^\La, e\in E_\La, \s^e\neq\s\}
\]
for $E_\La\coloneqq E\cap (\La\times\La)$,
and we define the \textit{cotangent bundle} of $S^\La$ by $T^*S^\La\coloneqq\Map(\Phi_\La,\R)$.
We call any element in the cotangent bundle a \textit{form}.
Suppose we are given a form $\omega\in \TS^\La$.  For any $e\in E_\La$, we 
define the function $\omega_e\in C(S^\La)$ by
\[
	\omega_e(\s)\coloneqq\omega((\s,\s^e))
\]
for any $\s\in S^\La$ such that $\s^e\neq\s$, and $\omega_e(\s)=0$ if $\s^e=\s$.
This gives a natural embedding
\begin{equation}\label{eq: 33}
	\TS^\La\hookrightarrow\prod_{e\in E_\La}C(S^\La), \qquad \omega\mapsto (\omega_e).
\end{equation}
The image of $\TS^\La$ in $\prod_{e\in E_\La}C(S^\La)$ corresponds to $(\omega_e)$ satisfying
$\omega_e(\s)=0$ if $\s^e=\s$,
and $\omega_e(\s)=\omega_{e'}(\s)$ if $\s^e=\s^{e'}$.
Since the system $\{ C(S^\La)\}_{\La\in\sI}$ form a projective system of $\R$-linear spaces
for the maps $\pi^\La\colon C(S^{\La'})\rightarrow C(S^\La)$ for $\La\subset \La'$ in $\sI$,
the product $\{\prod_{e\in E_\La} C(S^\La)\}_{\La\in\sI}$ 
also form a projective system.

\begin{lmm}\label{lem: 50}
	Consider $\La,\La'\in\sI$ such that $\La\subset\La'$.
	The projection $\pi^\La$ on the product
	induces an $\R$-linear homomorphism 
	$\pi^\La\colon\TS^{\La'}\rightarrow\TS^{\La}$.
\end{lmm}

\begin{proof}
	By definition of the embedding \eqref{eq: 33},
	we first
	prove that for any $e\in E_\La$, if $\omega_e(\s')=0$ for any $\s'\in S^{\La'}$ such that $\s'^e=\s'$,
	then $\pi^\La\omega_e(\s)=0$ for any $\s\in S^\La$ such that $\s^e=\s$.
	By calculation of the conditional expectation \eqref{eq: conditional expectation}, we have
	\begin{align*}
		\pi^\La\omega_e(\s)=\frac{1}{\mu(\s)}\sum_{\s'\in A_\s}\omega_e(\s').
	\end{align*}
	Since $\s=\s^e$, we have $\s'=\s'^e$ for any $\s'\in A_\s$.  
	This proves that $\pi^\La\omega_e(\s)=0$ as desired.
	Furthermore, if $\s^e=\s^{e'}$ for some $\s\in S^\La$ 
	and $e,e'\in E_\La\subset E_{\La'}$, then $\s'^{e}=\s'^{e'}$
	for any $\s'\in A_\s$, hence $\omega_{e}(\s')=\omega_{e'}(\s')$.
	This proves that $\pi^\La\omega_e(\s)=\pi^\La\omega_{e'}(\s)$ as desired.
\end{proof}

\begin{dfn}\label{def: cotangent}
	We define the space of \textit{co-local cotangent bundle}
	$\TS^X_\col$ by
	\[
		\TS^X_\col\coloneqq\varprojlim_{\sI} \TS^\La,
	\]
	where the limit is the projective limit with respect to $\pi^\La$.
\end{dfn}

We next define the space $C^1$ on $S^\La$.
Let
\[
	C^1(S^\La)\coloneqq \Map^\alt(\Phi_\La,\R),
\]
where
\[
	\Map^\alt(\Phi_\La,\R)\coloneqq\{\omega\in \TS^\La\mid \omega(\bar\vp)=-\omega(\vp)\}.
\]
Here $\bar\vp\coloneqq(t(\vp),o(\vp))$ for any $\vp=(o(\vp),t(\vp))\in\Phi_\La$.
We define the differential homomorphism $\partial_\La\colon C^0(S^\La)\rightarrow C^1(S^\La)$ by
\[
	\partial_\La f(\vp)\coloneqq f(t(\vp))-f(o(\vp))
\]
for any $f\in C(S^\La)$ and $\vp=(o(\vp),t(\vp))\in\Phi_\La$.
For any $\omega\in C^1(S^\La)$,
the image of $\omega$ with respect to the map \eqref{eq: 33} in $\prod_{e\in E}C(S^\La)$ 
consists of $(\omega_e)$ satisfying $\omega_e(\s)=0$ if $\s^e=\s$,
$\omega_e(\s)=\omega_{e'}(\s)$ if $\s^e=\s^{e'}$,
and $\omega_e(\s)=-\omega_{\bar e}(\s^e)$ if $\s^e\neq\s$.
The differential may be expressed as $\partial_\La f=(\nabla_e f)_{e\in E_\La}$,
where 
\[
	\nabla_e f(\s)\coloneqq f(\s^e)-f(\s)
\]
for any $e_\La\in E$ and $\s\in S^\La$.

Throughout the rest of this section, we assume the following.

\begin{dfn}\label{def: ordinary}
	In what follows, we say that a measure $\mu$ is \textit{ordinary}, if 
	for any $\La,\La'\in\sI$ such that $\La\subset\La'$
	and configurations $\s'\in S^{\La'}$, if we let $\s\coloneqq\pr_\La(\s')\in S^\La$ for the projection
	$\pr_\La\colon S^{\La'}\rightarrow S^\La$, then
	we have 
	\begin{equation}\label{eq: ordinary}
		\mu(\s^e)\mu(\s')=\mu(\s)\mu(\s'^e)
	\end{equation}
	for any $e\in E_{\La}$.  
\end{dfn}

\begin{lmm}\label{lem: ordinary}
	If $\mu$ is the product measure $\mu=\prod_{x\in X}\nu_x$ for some 
	family of probability measures $\{\nu_x\}_{x\in X}$ on $S$, then $\mu$ is ordinary.
\end{lmm}

\begin{proof}
	Let the notations be as in Definition \ref{def: ordinary}.
	If we let $\s=(\eta_x)\in S^\La$
	and $\s'=(\eta'_x)\in\La'$, then we have $\eta_x=\eta'_x$ for $x\in \La$, hence
	\begin{align*}
		\mu(\s)&=\prod_{x\in\La} \nu_x(\eta_x),&
		\mu(\s')&=\prod_{x\in\La} \nu_x(\eta_x) \times \prod_{x\in \La'\setminus\La}\nu_x(\eta'_x).
	\end{align*}
	Our assertion follows from the fact that since $e\in E_\La$, we have $\eta'^e_x=\eta'_x$ for any $x\in \La'\setminus\La$.
\end{proof}

If the probability measure $\mu$ is ordinary, then $\pi^\La$ preserves the local forms.

\begin{lmm}\label{lem: 505}
	If the probability measure $\mu$ is ordinary, then
	the projections $\pi^\La$ preserve the space of forms $C^1(S^\La)$ in $\TS^\La$.
	In other words, we have a homomorphism 
	\[	
		\pi^\La\colon C^1(S^{\La'})\rightarrow C^1(S^\La)
	\]
	for any $\La,\La'\in\sI$ such that $\La\subset\La'$.
\end{lmm}

\begin{proof}
	By Lemma \ref{lem: 50} and the definition of $C^1(S^\La)$ in Definition \ref{def: cotangent},
	it is sufficient to prove that if $(\omega_e)$ satisfies $\omega_e(\s')=-\omega_{\bar e}(\s'^e)$
	for any $\s'\in S^{\La'}$ such that $\s'^e\neq\s'$, 
	then $\pi^\La\omega_e(\s)=-\pi^\La\omega_{\bar e}(\s^e)$ for any $e\in E_\La$.
	Let $\s\in S^\La$ and $A_\s\coloneqq\{\s'\in S^{\La'}\mid \pr_\La(\s')=\s\}$, where
	$\pr_\La\colon S^{\La'}\rightarrow S^\La$ is the projection.  
	By calculation \eqref{eq: conditional expectation} of the conditional expectation, we have
	\begin{align*}
		\pi^\La\omega_e(\s)&
		=\frac{1}{\mu(\s)}\sum_{\s'\in A_\s}\omega_e(\s')\mu(\s'),&
		\pi^\La\omega_{\bar e}(\s^e)&
		=\frac{1}{\mu(\s^e)}\sum_{\s'\in A_{\s}}\omega_{\bar e}(\s'^e)\mu(\s'^e).
	\end{align*}
	Our assertion follows from \eqref{eq: ordinary}.
\end{proof}

By Lemma \ref{lem: 505},
we define the space of co-local forms $C^1_\col(S^X)$ as follows.

\begin{dfn}
	If the probability measure $\mu$ is ordinary, then
	we define the space of co-local forms by
	\[
		C^1_\col(S^X)\coloneqq\varprojlim_{\sI} C^1(S^\La),
	\]
	where the limit is the projective limit with respect to $\pi^\La$.
\end{dfn}

\begin{prp}\label{prop: differential}
	If the probability measure $\mu$ is ordinary, then
	the projection $\pi^\La$ is compatible with the differentials $\partial_\La$ on $C^0(S^\La)$.
	Hence the differential on each $C^0(S^\La)$ induces the differential
	\[
		\partial\colon C^0_\col(S^X)\rightarrow C^1_\col(S^X)
	\]
	on the space of co-local functions.
\end{prp}

\begin{proof}
	It is sufficient to prove that for any $\La,\La'\in\sI$ such that $\La\subset\La'$,
	we have $\pi^\La(\nabla_e f)=\nabla_e(\pi^\La f)(\s)$
	for any $e\in E_\La$ and $f\in C^0(S^{\La'})$. 
	By definition of the differential and calculation of the conditional expectation \eqref{eq: conditional expectation},
	we have
	\begin{align*}
		\pi^\La(\nabla_e f)(\s)&=\frac{1}{\mu(\s)} \sum_{\s'\in A_\s}\nabla_e f(\s')\mu(\s')
		=\frac{1}{\mu(\s)} \sum_{\s'\in A_\s} f(\s'^e) \mu(\s') -
		\frac{1}{\mu(\s)} \sum_{\s'\in A_\s} f(\s') \mu(\s'),
	\end{align*}
	where $A_\s=\{\s'\in S^{\La'}\mid\pr_\La(\s')=\s\}$ as before.
	Similarly, we have
	\begin{align*}
		\nabla_e(\pi^\La f)(\s)&=(\pi^\La f)(\s^e)- (\pi^\La f)(\s)
		=\frac{1}{\mu(\s^e)}\sum_{\s'\in A_{\s}} f(\s'^e)\mu(\s'^e)
		-\frac{1}{\mu(\s)}\sum_{\s'\in A_\s} f(\s')\mu(\s').
	\end{align*}
	Our assertion now follows from the equality \eqref{eq: ordinary} satisfied by ordinary probability measures.
\end{proof}

\begin{dfn}\label{def: colocal}
	If the probability measure $\mu$ is ordinary,
	then we let
	\[
		H^0_\col(S^X)\coloneqq \Ker\bigl(\partial\colon C^0_\col(S^X)\rightarrow C^1_\col(S^X)\bigr)
	\]
	for the differential $\partial\colon C^0_\col(S^X)\rightarrow C^1_\col(S^X)$
	of Proposition \ref{prop: differential}.
\end{dfn}

Next we define the notion of closed local and co-local forms.
For any $\Lambda\subset X$, we define a path in $S^\La$ to be
a sequence of transitions $\vec\gamma=(\vp^1,\ldots,\vp^N)$ in $\Phi_\La$
such that $t(\vp^i)=o(\vp^{i+1})$ for $0<i<N$.  We say that the path is 
\textit{closed}, if $t(\vp^N)=o(\vp^1)$. 
For a form $\omega\in C^1(S^\La)$, we define the integration with respect
to a path $\vec\gamma=(\vp^1,\ldots,\vp^N)$ by
\[
	\int_{\vec\gamma}\omega\coloneqq\sum_{i=1}^N\omega(\vp^i).
\]
As in \cite[Definition 2.14]{a} we define the closed forms as follows.

\begin{dfn}
	We say that a form $\omega\in C^1(S^\La)$ is \textit{closed}, if 
	\[
		\int_{\vec\gamma}\omega=0
	\]
	for any closed path $\vec\gamma$ in $S^\La$.
\end{dfn}

We denote by $Z^1(S^\La)$ the $\R$-linear space of closed forms in $C^1(S^\La)$.
The path $\vec\gamma=(\vp^1,\ldots,\vp^N)$ may be written as $\vec\gamma=(\s^0,\ldots,\s^N)$,
where $\vp^i=(\s^{i-1},\s^i)$ for $1\leq i\leq N$.
By definition of a transition, there exists $e_1,\ldots,e_N\in E_\La$ 
such that $\s^{i}=(\s^{i-1})^{e_{i}}$ for $1\leq i\leq N$.
The path is closed if $\s^0=\s^N$.   If $\omega=(\omega_e)$ is a closed form in $C^1(S^\La)$, then we have
\[
	\sum_{i=1}^N \omega_{e_i}(\s^{i-1})=0.
\]

\begin{lmm}
	If the probability measure $\mu$ is ordinary, then
	the projection $\pi^\La$ induces an $\R$-linear homomorphism
	\[
		\pi^\La\colon Z^1(S^{\La'})\rightarrow Z^1(S^\La)
	\]
	for any $\La,\La'\in\sI$ such that $\La\subset\La'$.
\end{lmm}

\begin{proof}
	Let $\omega=(\omega_e)$ be a closed form in $Z^1(S^{\La'})$, 
	and consider the form $\pi^\La\omega$, which is an element in 
	$C^1(S^\La)$ by Lemma \ref{lem: 505}.
	Note that we have $\pi^\La\omega=(\pi^\La\omega_e)$.
	Consider a closed path $\vec\gamma=(\s^0,\ldots,\s^N)$ in $S^\La$
	such that $\s^{i}=(\s^{i-1})^{e_{i}}=\cdots =(\s^0)^{e_1\cdots e_i}$ for some $e_1,\ldots,e_N\in E_\La$.
	Note that for any $\tilde\eta^0\in A_{\s_0}\subset S^{\La'}$, if we let 
	$\tilde\s^{i}=(\tilde\s^0)^{e_1\cdots e_i}$ for $1\leq i\leq N$,
	then $\vec\gamma'\coloneqq(\tilde\s^0,\ldots,\tilde\s^N)$ is a closed path in $S^{\La'}$.
	Then we have
	\begin{align*}
		\int_{\vec\gamma}\pi^\La\omega=\sum_{i=1}^N\pi^\La\omega_{e_i}(\s^{i-1})
		&=\sum_{i=1}^N\frac{1}{\mu(\s^{i-1})}\sum_{\tilde\s^0\in A_{\s_0}}\omega_{e_i}(\tilde\s^{i-1})\mu(\tilde\s^{i-1})\\
		&=\sum_{\tilde\s^0\in A_{\s^0}}\frac{\mu(\tilde\s^0)}{\mu(\s^0)}\left(\sum_{i=1}^N\omega_{e_i}(\tilde\s^{i-1})\right)=0,
	\end{align*}
	where the third equality follows from \eqref{eq: ordinary},
	and the last equality follows from the fact that $\omega$ is closed in $S^{\La'}$.
	This gives our assertion.
\end{proof}

We now give the definition of closed co-local forms.

\begin{dfn}
	Assume that the probability measure $\mu$ is ordinary.
	We define the space of closed co-local forms $Z^1_\col(S^X)$ by
	\begin{align*}
		Z^1_\col(S^X)&\coloneqq\varprojlim_{\La\in\sI}Z^1_\col(S^\La),
	\end{align*}
	where the limit is the projective limit with respect to $\pi^\La$.
\end{dfn}

We next consider the integration of closed forms.
For any $\La\in\sI$, by \cite[Lemma 2.14, Lemma 2.15]{a},  we have 
an exact sequence
\begin{equation}\label{eq: SES3}
	\xymatrix{
	0\ar[r]& \Ker\partial_\La\ar[r]&C^0(S^\La)\ar[r]^{\partial_\La}& Z^1(S^\La)\ar[r]&0.
	}
\end{equation}
Assume that the probability measure $\mu$ is ordinary.
By Proposition \ref{prop: differential}, 
since the differential $\partial_\La$ is compatible with the projection $\pi^\La$,
we see that $\{\Ker\partial_\La\}_{\La\in\sI}$ also form a projective system for the projection $\pi^\La$.
The Mittag-Leffler condition for projective systems is given as follows (see for example \cite[(13.1.2)]{EGAIII}).

\begin{dfn}
	A projective system $\{M_\La\}_{\La\in\sI}$ satisfies the \textit{Mittag-Leffler condition}, if for any $\La\in \sI$, if we let
	$N_\La\coloneqq\bigcap_{\La'\in\sI, \La\subset \La'}  \Im(\pi^{\La}\colon M_{\La'}\rightarrow M_{\La})$, 
	then there exists $\La''\in\sI$ satisfying $\La\subset \La''$
	such that $N_\La= \Im(\pi^{\La}\colon M_{\La''}\rightarrow M_{\La})$.  
	In other words, the image of $\pi^{\La}$ stabilizes for $\La''$ sufficiently large.
\end{dfn}

By \cite[Proposition (13.2.2)]{EGAIII},
if the projective system $\{\Ker\partial_\La\}_{\La\in\sI}$ 
satisfies the Mittag-Leffler condition, then the projective limit
\[
	\xymatrix{
	0\ar[r]& \varprojlim_{\sI}\Ker\partial_\La\ar[r]&
	\varprojlim_{\sI} C^0(S^\La)\ar[r]^{\partial_\La}&\varprojlim_{\sI}
	Z^1(S^\La)\ar[r]&0
	}
\]
of \eqref{eq: SES3} is exact.  We will use this fact to prove the following.

\begin{prp}\label{prop: important}
	If the probability measure $\mu$ on $S^X$ is ordinary, then
	the projective system $\{\Ker\partial_\La\}_{\La\in\sI}$ satisfies
	the Mittag-Leffler condition.  This implies that the sequence
	\begin{equation}\label{eq: SES}
		\xymatrix{
		0\ar[r]& H^0_\col(S^X)\ar[r]&C^0_\col(S^X)\ar[r]^{\partial}& Z^1_\col(S^X)\ar[r]&0
	}
	\end{equation}
	is exact, where $H^0_\col(S^X)=\varprojlim_{\sI}\Ker\partial_\La$ as in Definition \ref{def: colocal}.
\end{prp}

\begin{proof}
	We first note that by \cite[Remark 2.28]{a}, any function $f\in C(S^\La)$
	satisfies $\partial_\La f=0$ if and only if $f$ is constant on each of the connected components of
	the graph $(S^\La,\Phi_\La)$.  
	Since $\La\in\sI$ is finite, the configuration space $S^\La$ has only a finite number of connected
	components.  This shows that $\Ker\partial_\La$ is finite dimensional.
	This ensures that $\{\Ker\partial_\La\}_{\La\in\sI}$ satisfies the Mittag-Leffler condition, since
	any descending sequence of linear subspaces of a finite dimensional linear space is stable.
	Our assertion now follows from \cite[Proposition (13.2.2)]{EGAIII}.
\end{proof}


%
%
%
\section{Decomposition of Varadhan Type}\label{sec: Varadhan}
%
%
%

In this section, we consider a group $G$ and an action of $G$ on the locale $X$,
and prove the decomposition theorem of Varadhan type for closed co-local forms.
Before going into the details, we first fix notations concerning action of $G$.
An action of $G$ on $X$ gives a bijection $\sigma\colon X\rightarrow X$ for any $\sigma\in G$.
We define the action of $G$ on $S^X$ given by mapping $\s=(\eta_x)_{x\in X}\in\prod_{x\in X}S$ to 
$\sigma(\s)\coloneqq(\eta_{\sigma^{-1}(x)})_{x\in X}$ for any $\sigma\in G$.
Then $G$ induces an action on $C(S^X)$ given for any $\sigma\in G$ by
\[
	\sigma(f)(\s)=f(\sigma^{-1}(\s))
\]
for any $\s\in S^X$.  

For any subset $\La\subset X$, the element $\sigma\in G$ induces a bijection 
$\sigma\colon\La\cong\sigma(\La)$.  Hence on $S^\La$, this induces a bijection
\begin{equation}\label{eq: s}
	\sigma\colon S^{\La}\cong S^{\sigma(\La)}.
\end{equation}
Note that for any $\s\in S^{\La}$, we have $\sigma(\s)\in S^{\sigma(\La)}$.
In terms of components, if we let $\s=(\eta_x)_{x\in\La}\in S^{\La}$, 
then we have $\sigma(\s)=(\eta_{\sigma^{-1}(x)})_{x\in\sigma(\La)}\in S^{\sigma(\La)}$.
Hence \eqref{eq: s} induces a bijection
\begin{equation}\label{eq: function}
	\sigma\colon C\bigl(S^\La\bigr)\cong C\bigl(S^{\sigma(\La)}\bigr),
\end{equation}
which maps any function $f\in C\bigl(S^{\La}\bigr)$ to the function $\sigma(f)\in C\bigl(S^{\sigma(\La)}\bigr)$.
For any $\La\in\sI$, we have $\sigma(\La)\in\sI$ for any $\sigma\in G$.
Clearly, the action of $G$ is compatible with the natural inclusion $C(S^\La)\hookrightarrow C(S^X)$, hence
we have an action of $G$ on $C_\loc(S^X)$.

The action of the group $G$ on $X$ defines a map of graphs 
$\sigma\colon (S^\La,\Phi_\La)\rightarrow (S^{\sigma(\La)},\Phi_{\sigma(\La)})$,
which induces a natural actions on $C(\Phi_\La)$ by $\sigma(\omega)(\vp)=\omega(\sigma^{-1}(\vp))$.
We have the following.
\begin{lmm}\label{lem: differential}
	For any $\La\in\sI$,
	the action of $G$ is compatible with the 
	differential 
	\[
		\partial_\La\colon C^0(S^\La)\rightarrow C^1(S^\La).
	\]
\end{lmm}

\begin{proof}
	For $f\in C(S^\La)$, we have
	\begin{align*}
		\sigma(\partial f)(\vp)&=\partial f(\sigma^{-1}(\vp))=f(\sigma^{-1}(t(\vp)))-f(\sigma^{-1}(o(\vp)))\\
		&=\sigma(f)(t(\vp))-\sigma(f)(o(\vp))=(\partial \sigma(f))(\vp).
	\end{align*}
	This shows that $\partial_\La$ is compatible with the action of the group $G$ as desired.
\end{proof}

Note that if we let $\cF_X$ be the Borel $\sigma$-algebra for $S^X$, then we have
$\sigma(A)\in\cF_X$ for any $A\in\cF_X$ and $\sigma\in G$.  We define a $G$-invariant
probability measure on $S^X$ as follows.
\begin{dfn}
	Suppose $\mu$ is a probability measure on $S^X$.
	We say that $\mu$ is \textit{invariant with respect to the action of $G$},
	if we have 
	\[
		\mu(\sigma(A))=\mu(A)
	\]
	for any $A\in\cF_X$.
\end{dfn}

If $\mu=\prod_{x\in X}\nu_x$ for some family of probability measures $\{\nu_x\}$ on $S$,
and if $\nu_{\sigma(x)}=\nu_x$ for any $\sigma\in G$, then the probability measure
$\mu$ is invariant with respect to the action of $G$.
In particular, the product measure $\mu=\nu^{\otimes X}$ for a measure $\nu$ on $S$ is 
invariant with respect to the action of $G$.

For the rest of this section, we assume that $\mu=\nu^{\otimes X}$ for some probability measure $\nu$ on $S$ 
which is supported on $S$.  
 Then $\mu$ is invariant with respect to the action of $G$.
For any $f\in C(S^X)$ and $\La\in\sI$, if we let $\pi^\La f\in C(S^\La)$ be the conditional expectation,
then we have $\sigma(\pi^\La f)\in C(S^{\sigma(\La)})$ for any $\sigma\in G$.
Moreover, since the action of $\sigma\in G$ is compatible with
 all of the structures appearing in the conditional expectation, we have
 \[
 	\sigma(E_\mu[f|\cF_\La])=E_{\sigma(\mu)}[\sigma(f)|\cF_{\sigma(\La)}].
 \]
 Our condition that $\mu$ is invariant with respect to the action of $G$ gives the equality
$\sigma(\pi^\La f)=\pi^{\sigma(\La)}(\sigma(f))$.  This shows that the projections
$\pi^\La$ are compatible with the action of $G$,
hence we have an action of $G$ on the co-local functions
$C_\col(S^X)=\varprojlim_{\sI}C(S^\La)$. 
The action \eqref{eq: function} gives a mapping
\[
	\sigma\colon C_{\La}(S^X)\rightarrow C_{\sigma(\La)}(S^X)
\]
on $C_{\La}(S^X)$.
By Proposition \ref{prop: expansion},
we have a unique expansion
\[
	f = \sum_{\La\in\sI} f_{\La}
\]
for any $f\in C_{\col}(S^X)$, where $f_{\La}\in C_{\La}(S^X)$ for any $\La\in\sI$.
Then we have
\[
	\sigma(f) =  \sum_{\La\in\sI}\sigma(f_{\La}).
\]
Since $f_{\La}\in C_{\La}(S^X)$, we have $\sigma(f_{\La})\in C_{\sigma(\La)}(S^X)$,
hence the uniqueness of expansion gives
\[
	\sigma(f_{\La})=\sigma(f)_{\sigma(\La)}
\]
for any $\La\in\sI$.   This proves the following assertion.

\begin{lmm}
	The action of the group $G$ on the locale $X$ gives a natural action of $G$
	on the space of uniformly co-local functions $C^0_\unif(S^X)$.
\end{lmm}

The probability measure $\mu=\nu^{\otimes X}$ is invariant by the action of $G$.
In addition, by Lemma \ref{lem: ordinary}, the measure $\mu$ is ordinary.
By Proposition \ref{prop: important}, we have an exact sequence
\[
	\xymatrix{
	0\ar[r]& H^0_\col(S^X)\ar[r]&C^0_\col(S^X)\ar[r]^{\partial}& Z^1_\col(S^X)\ar[r]&0}.
\]
Considering the long exact sequence associated to the group cohomology of $G$, we have an
exact sequence
\begin{equation}\label{eq: LES}
	\cdots\rightarrow C^0_\col(S^X)^G\rightarrow Z^1_\col(S^X)^G\xrightarrow\delta
	H^1_\col(G, H^0_\col(S^X))\rightarrow H^1_\col(G,C^0_\col(S^X))\rightarrow\cdots,
\end{equation}
where $C^0_\col(S^X)^G$ and $Z^1_\col(S^X)^G$ denote the $G$-invariant subgroups of
$C^0_\col(S^X)$ and $Z^1_\col(S^X)$, and
$H^1(G,H^0_\col(S^X))$ and $H^1(G,C^0_\col(S^X))$ are group cohomology of $G$ with coefficients in $H^0_\col(S^X)$ and $C^0_\col(S^X)$.
This gives an injective homomorphism 
\begin{equation}\label{eq: boundary}\xymatrix{
	Z^1_\col(S^X)^G/\partial\bigl(C^0_\col(S^X)^G\bigr)\,\ar@{^(->}[r]& \,H^1(G, H^0_\col(S^X)).
}\end{equation}
The boundary morphism $\delta$ is given explicitly by mapping any $\omega\in Z^1_\col(S^X)^G$
to the cocycle $\rho$ given by
$\rho(\sigma)=(1-\sigma)\theta$ for any $\sigma\in G$, 
where $\theta\in C^0_\col(S^X)$ is a function satisfying
$\partial\theta=\omega$.

Next, suppose the group $G$ is torsion free and the action of $G$ on $X$ is free.
Then the action of $G$ 
on the set $\sI\setminus\{\emptyset\}$ is free.  We denote by $\sI_0$ a set of representatives of 
the equivalence classes of
$\sI\setminus\{\emptyset\}$ with respect to the action of $G$.  This implies that for any nonempty
$\La\subset X$, there exists a unique $\tau\in G$ such that $\tau^{-1}(\La)\in\sI_0$.
For any $f\in C^0_\col(S^X)$, the canonical decomposition
\[
	f=\sum_{\La\in\sI} f_\La,\qquad f_\La\in C_\La(S^X)
\]
with $f_\emptyset=0$ may be written as
\[
	f=\sum_{\tau\in G}\sum_{\La_0\in\sI_0}f_{\tau(\La_0)}.
\]
For any $\tau\in G$, if we let
\begin{equation}\label{eq: fs}
	f_\tau\coloneqq\sum_{\La_0\in\sI_0}f_{\tau(\La_0)},
\end{equation}
then we have
\[
	f=\sum_{\tau\in G} f_\tau.
\]

\begin{thm}\label{thm: main}
	Suppose $\mu=\nu^{\otimes X}$, and that the group
	$G$ is torsion free and the action of $G$ on $X$ is free.  
	Then the boundary morphism \eqref{eq: boundary} induces an isomorphism
	\begin{equation}\label{eq: isom}
		Z^1_\col(S^X)^G/\partial\bigl(C^0_\col(S^X)^G\bigr)\cong H^1(G, H^0_\col(S^X)).
	\end{equation}
	In particular, the choice of $\sI_0$ gives a splitting of the boundary morphism,
	hence a decomposition
	\begin{equation}\label{eq: decomposition}
		Z^1_\col(S^X)^G\cong \partial\bigl(C^0_\col(S^X)^G\bigr) \oplus H^1(G, H^0_\col(S^X)).
	\end{equation}
\end{thm}

\begin{proof}
	It is sufficient to construct a section of the boundary morphism
	\begin{equation}\label{eq: boundary2}
		\delta\colon Z^1_\col(S^X)^G \rightarrow  H^1(G, H^0_\col(S^X)).
	\end{equation}
	Let $\rho\in Z^1(G, H^0_\col(S^X))$ be a group 
	cocycle representing a class in $H^1(G, H^0_\col(S^X))$.
	Then $\rho$ is a map from $G$ to $H^0_\col(S^X)=\Ker\partial\subset C^0_\col(S^X)$ 
	satisfying $\rho(\sigma\tau)=\sigma\rho(\tau)+\rho(\sigma)$ for any $\sigma,\tau\in G$.
	We let
	\[
		\theta_\rho\coloneqq \sum_{\tau\in G}\rho(\tau)_\tau\in C^0_\col(S^X),
	\]
	where $\rho(\tau)_\tau$ is the function $f_\tau$ in \eqref{eq: fs} for $f=\rho(\tau)$.
	Then we have
	\begin{align*}
		\sigma(\theta_\rho)&=\sum_{\tau\in G}\sigma(\rho(\tau)_\tau)
		=\sum_{\tau\in G}\sigma(\rho(\tau))_{\sigma\tau}
		=\sum_{\tau\in G}(\rho(\sigma\tau)_{\sigma\tau}-\rho(\sigma)_{\sigma\tau})\\
		&=\sum_{\tau\in G}(\rho(\tau)_{\tau}-\rho(\sigma)_{\tau})
		=\theta_\rho-\rho(\sigma).
	\end{align*}
	This shows that $(1-\sigma)\theta_\rho=\rho(\sigma)$.
	Since $\rho(\sigma)\in H^0_\col(S^X)=\Ker\partial$, the compatibility of $\partial$ with the
	action of $G$ gives
	\[
			(1-\sigma)\partial \theta_\rho=\partial ((1-\sigma)\theta_\rho)=\partial \rho(\sigma)=0
	\]
	for any $\sigma\in G$,
	which shows that $\omega_\rho\coloneqq\partial\theta_\rho\in Z^1_\col(S^X)^G$.
	By the definition of the boundary morphism, we see that $\delta(\omega_\rho)=\rho$
	as an element in $H^1(G, H^0_\col(S^X))$.  This shows that the map $\rho\mapsto\omega_\rho$
	gives a section of \eqref{eq: boundary2} as desired.
	This proves that \eqref{eq: isom} is in fact an isomorphism and gives the decomposition
	\eqref{eq: decomposition} for our choice of $\sI_0$.
\end{proof}

We next consider the property that an interaction is irreducibly quantified,
which was originally proposed in \cite[Definition 2.22]{a}.

\begin{dfn}
	We say that the interaction $(S,\phi)$ is \textit{irreducibly quantified},
	if for any locale $(X,E)$ whose set of vertices is \textit{finite}, 
	the associated configuration space with transition structure
	$S^X$ satisfies the following property. For any $\s,\s'\in S^X$, if $\xi_X(\s)=\xi_X(\s')$
	for any conserved quantity $\xi\in\C{S}$, then there exists a path $\vec\gamma$ from
	$\s$ to $\s'$ in $S^X$.  Here, $\xi_X(\s)\coloneqq\sum_{x\in X}\xi(\eta_x)$ for any $\s=(\eta_x)\in S^X$,
	which is a finite sum since we have assumed that $X$ is finite.
\end{dfn}	

Since we have assumed that $S$ is finite, the dimension $c_\phi\coloneqq\dim_\R\C{S}$ is finite.
The following result follows from \cite[Lemma 2.28]{a} and the definition of irreducibly quantified
interactions. We will include the proof for the sake of completeness.

\begin{lmm}\label{lem: review}
	Let $(X,E)$ be a finite locale, and suppose $(S,\phi)$ is an interaction which is irreducibly quantified.
	Suppose that $f\in C(S^X)$ satisfies $\partial_X f=0$.
	Then we have $f(\s)=f(\s')$
	for any $\s,\s'\in S^X$ satisfying $\xi_X(\s)=\xi_X(\s')$ for any conserved quantity $\xi\in\C{S}$.
\end{lmm}

\begin{proof}
	Let $\s,\s'\in S^X$, and suppose that $\xi_X(\s)=\xi_X(\s')$
	for any conserved quantity $\xi\in\C{S}$.
	Since the pair $(S,\phi)$ is irreducibly quantified, 
	there exists a path $\vec\gamma$ from $\s$ to $\s'$ in $S^X$.
	Write $\vec\gamma=(\varphi^1,\ldots,\varphi^N)$
	with transitions $\varphi^i=(\s^{i-1},\s^i)\in\Phi$ for $i=1,\ldots,N$.
	Since $\partial_X f=0$, we have
	\[
		f(\s^i)-f(\s^{i-1})=\partial_X f(\vp^i)=0
	\]
	for any $i=1,\ldots,N$.  This proves that $f(\s)=f(\s^0)=f(\s^N)=f(\s')$ as desired.
\end{proof}

Using Lemma \ref{lem: review}, we may prove the following.

\begin{lmm}\label{lem: invariant}
	Suppose the interaction $\phi$ is irreducibly quantified.  Then
	we have 
	\[
		H^0_\col(S^X)\subset C^0_\col(S^X)^G.
	\]
\end{lmm}

\begin{proof}
	Let $f=(f^\La)\in H^0_\col(S^X)$. 
	Note that for any $\sigma\in G$, we have $\sigma(f)=(\sigma(f^\La))$, 
	where $\sigma(f^\La)\in C(S^{\sigma(\La)})$.
	It is sufficient to prove that $\sigma(f^\La)=f^{\sigma(\La)}$.
	By definition of $H^0_\col(S^X)$ given in Definition \ref{def: colocal}, 
	we have $f^\La\in \Ker\partial_\La$.
	Let $\sigma\in G$, and take $\La'\in\sI$ sufficiently large so that $\La,\sigma(\La)\subset\La'$.
	We fix an arbitrary bijection $\upsilon\colon\La'\rightarrow\La'$ 
	satisfying $\upsilon(x)=\sigma^{-1}(x)$ if $x\in\sigma(\La)$.
	 We define a map $T_{\upsilon}\colon S^{\La'}\rightarrow S^{\La'}$ by 
	\[
		(T_{\upsilon}(\s))_x \coloneqq\eta_{\upsilon(x)}.
	\]
	By construction, $\xi_{\La'}(\s)=\xi_{\La'}(T_{\upsilon}(\s))$.
	Since $\partial_{\La'}f^{\La'}=0$, by Lemma \ref{lem: review}
	and the fact that the interaction $\phi$ is irreducibly quantified,
	we have
	$
		f^{\La'}(\s)=f^{\La'}(T_{\upsilon}(\s))
	$
	for any $\s\in S^{\La'}$.   This shows that
	\[
		\sigma(f^\La)=\sigma(\pi^\La f^{\La'})=\pi^{\sigma(\La)}(f^{\La'}\circ
		T_{\upsilon})=\pi^{\sigma(\La)}f^{\La'}=f^{\sigma(\La)},
	\]
	which gives the desired result.
\end{proof}

Using this fact, we may deduce the following from Theorem \ref{thm: main}.

\begin{crl}\label{crl: main}
	Assume the conditions of Theorem \ref{thm: main} and suppose in addition that the interaction
	$\phi$ is irreducibly quantified.  Then we have an isomorphism
	\[
		Z^1_\col(S^X)^G/\partial\bigl(C^0_\col(S^X)^G\bigr)\cong
		\Hom_\bbZ(G,H^0_\col(S^X)).
	\]
	In particular, if the maximal abelian quotient of
	$G$ is finitely generated of rank $d$, then we have a decomposition
	\[
		Z^1_\col(S^X)^G\cong\partial\bigl(C^0_\col(S^X)^G\bigr)\oplus\bigoplus_{j=1}^dH^0_\col(S^X).
	\]
\end{crl}

\begin{proof}
	By Lemma \ref{lem: invariant}, the action of $G$ on $H^0_\col(S^X)$ is trivial.  Hence we have
	\[
		H^1(G,H^0_\col(S^X))\cong \Hom_\bbZ(G,H^0_\col(S^X)),
	\]
	which combined with Theorem \ref{thm: main} gives the first isomorphism.
	If we fix a generator $\tau_1,\ldots,\tau_d$ of the maximal abelian quotient of $G$, 
	then an element of 
	$\Hom_\bbZ(G,H^0_\col(S^X))$ is determined by the images of $\tau_i$ in $H^0_\col(S^X)$.
	This fact and the decomposition \eqref{eq: decomposition}
	of Theorem \ref{thm: main} gives the second isomorphism as desired.
\end{proof}

%
%
%
\section{The $L^2$-Case: A Conjecture}\label{sec: L2}
%
%
%

The decomposition theorem for Varadhan that is necessary for proving the hydrodynamic limit
is a decomposition for $L^2$-forms.
In this section, we give the definition of $L^2$-forms and formulate a conjecture
concerning Varadhan's decomposition in this case.
Let $\mu$ be a probability measure on $S^X$ supported on $S^\La$ for any $\La\in\sI$.
We let
\[
	\dabs{f}_\mu\coloneqq\pair{f,f}^{1/2}_\mu=E_\mu[f^2]^{1/2}
\]
for any $f\in C(S^X)$.

\begin{dfn}
	We define the  \textit{space of $L^2$-functions} $L^2(\mu)$ to be the quotient space
	\[
		L^2(\mu)\coloneqq\{f\in C(S^X)\mid \dabs{f}_\mu<\infty\}/
		\{ f\in C(S^X)\mid \dabs{f}_\mu=0\}.
	\]
\end{dfn}

By standard facts concerning $L^2$-spaces, $L^2(\mu)$ is known to be
a Hilbert space for the inner product $\pair{f,g}_\mu=E_\mu[fg]$.
In particular,  $L^2(\mu)$ is complete for the topology given by 
the norm $\dabs{\cdot}_\mu$.  

\begin{rem}
	It is well-known that  if $f\in L^2(\mu)$, then $f$ is integrable for $\mu$.
	This may be seen from the fact that since $\abs{f}\leq(f^2+1)/2$, we have
	\[
		E_\mu[\abs{f}]\leq E_\mu[(f^2+1)/2] =(E_\mu[f^2]+1)/2\leq\infty.
	\]
\end{rem}

For any $f\in L^2(\mu)$, the system $(\pi^\La f)_{\La\in\sI}$ defines a co-local function
in $C_\col(S^X)$,
hence we have a natural homomorphism
\begin{equation}\label{eq: nh}
	L^2(\mu)\rightarrow C_\col(S^X), \qquad f\mapsto (\pi^\La f)_{\La\in\sI}.
\end{equation}
Furthermore, if $f\in C(S^\La)$ for $\La\in \sI$, then
the function $f$ viewed as an element in $C(S^X)$ satisfies $\dabs{f}_\mu<\infty$, hence $f$
defines an element in $L^2(\mu)$.  This gives a natural homomorphism
\[
	C_\loc(S^X)\rightarrow L^2(\mu),
\]
which is injective since  the composite $C_\loc(S^X)\rightarrow L^2(\mu)\rightarrow C_\col(S^X)$
gives the natural inclusion.
Note that since $\pi^\La$ is an orthogonal projection (see Lemma \ref{lem: orthogonal}),
we have 
\begin{equation}\label{eq: bound}
	\dabs{\pi^\La f}_\mu =\pair{\pi^\La f,\pi^\La f}_\mu\leq\pair{f,f}_\mu=\dabs{f}_\mu
\end{equation}
for any $f\in L^2(\mu)$.
The following is the \textit{Martingale Convergence Theorem} for our case.

\begin{thm}\label{thm: convergence}
	We let
	\[
		C_{L^2}(S^X)\coloneqq\bigl\{ (f^\La)\in C_\col(S^X)\mid \sup_{\La\in\sI}\dabs{f^\La}_\mu<\infty\bigr\}.
	\]
	Consider a family of sets $\{\La_n\}_{n\in\bbN}$ in $\sI$ such that $\La_n\subset\La_{n+1}$ and 	
	$X=\bigcup_{n\in\bbN}\La_n$.
	For any $f\in C_{L^2}(S^X)$, if we let $f_n\coloneqq f^{\La_n}$ for any $n\in\bbN$, then 
	the sequence of functions $(f_n)_{n\in\bbN}$ converges strongly
	to a function $f_\infty$ in $L^2(\mu)$.
\end{thm}

\begin{proof}
	For $m>n$, since $f_n=\pi^{\La_n} f_m$, 
	we see from \eqref{eq: bound} that $\dabs{f_m}^2_\mu\geq \dabs{f_n}^2_\mu$,
	hence $\dabs{f_n}^2_\mu$ is monotonously increasing for $n\geq0$.
	Let
	\[
		M\coloneqq\sup_{n\in\bbN}\dabs{f_n}_\mu\leq \sup_{\La\in\sI}\dabs{f^\La}_\mu<\infty.
	\]
	By Lemma \ref{lem: orthogonal}, the projection $\pi^{\La_n}$ is orthogonal. This
	shows that we have
	\[
		\dabs{f_m}^2_\mu
		=\dabs{f_n}^2_\mu+\dabs{f_m-f_n}^2_\mu
	\]
	for any $m>n$, which gives the equality 
	$
		\dabs{f_m-f_n}^2_\mu=\dabs{f_m}^2_\mu-\dabs{f_n}^2_\mu.
	$
	This shows that
	\[
		\lim_{m,n\rightarrow\infty}\dabs{f_m-f_n}_\mu=\lim_{m\rightarrow\infty}\dabs{f_m}_\mu-
		\lim_{n\rightarrow\infty}\dabs{f_n}_\mu=M-M=0,
	\]
	hence $(f_n)_{n\in\bbN}$ is a Cauchy sequence in $L^2(\mu)$.
	Since $L^2(\mu)$ is complete for the $L^2$-norm,
	the sequence $(f_n)_{n\in\bbN}$ converges strongly to an element $f_\infty$ in $L^2(\mu)$ as desired.
\end{proof}

\begin{rem}
	Since $\sI$ is a directed set, the system of functions
	$(f^\La)\in C_{L^2}(S^X)$ may be interpreted as a \textit{net}
	with set of indices $\sI$.
	As a generalization of Theorem \ref{thm: convergence}, we may prove
	that for any $(f^\La)\in C_{L^2}(S^X)$, we have
	\[
		\lim_{\La\in\sI} f^\La=f_\infty\in L^2(\mu),
	\]
	where the limit is the convergence in terms of nets.
\end{rem}

Theorem \ref{thm: convergence} gives the following corollary.

\begin{crl}
	The space $C_{L^2}(S^X)$ of Theorem \ref{thm: convergence} coincides
	with the image of $L^2(\mu)$ in $C_\col(S^X)$ with respect to the homomorphism \eqref{eq: nh}.  	
	Moreover, the space of local functions
	$C_\loc(S^X)$ is dense in $C_{L^2}(S^X)$ for the topology defined by the
	norm
	\[
		\dabs{(f^\La)}_\mu\coloneqq\sup_{\La\in\sI}\dabs{f^\La}_\mu
	\]
	on $C_{L^2}(S^X)$.
\end{crl}

\begin{proof}
	Suppose $f\in L^2(\mu)$.  Then for any $\La\in\sI$, by \eqref{eq: bound}, we have
	\[
		\dabs{f^\La}_\mu\leq\dabs{f}_\mu,
	\]
	hence we see that $(\pi^\La f)$ satisfies $\sup_{\La\in\sI}\dabs{f^\La}_\mu \leq \dabs{f}_\mu<\infty$.
	This shows that $(\pi^\La f)\in C_{L^2}(S^X)$ as desired.
	On the other hand, suppose $(f^\La)\in  C_{L^2}(S^X)$.
	We fix a family $\{\La_n\}_{n\in\bbN}$ of sets in $\sI$ such that $\La_n\subset\La_{n+1}$ and $X=\bigcup_{n\in\bbN}\La_n$, and let $f_n\coloneqq f^{\La_n}$ for any $n\in\bbN$.
	Then Theorem \ref{thm: convergence} shows that $(f_n)$ converges strongly to a function $f_\infty$
	in $L^2(\mu)$.
	Fix an $n\in\bbN$ and let $g_n\in C(S^{\La_n})$.  
	Then for any $m>n$, the orthogonal property of the conditional 
	expectation in Lemma \ref{lem: orthogonal} implies that $\pair{f_m-f_n,g_n}_\mu=0$,
	hence we have
	\[
		\pair{f_m,g_n}_\mu=\pair{f_n+(f_m-f_n),g_n}_\mu=\pair{f_n,g_n}_\mu.
	\]
	This shows that 
	\[
		\pair{f_\infty,g_n}_\mu=\lim_{m\rightarrow\infty}\pair{f_m,g_n}_\mu=\pair{f_n,g_n}_\mu.
	\]
	Furthermore, the definition of conditional expectation shows that we have
	\[
		\pair{\pi^{\La_n}f_\infty,g_n}_\mu=\pair{f_\infty,g_n}_\mu.
	\]
	This proves that $\pi^{\La_n}f_\infty=f_n$ in $C(S^{\La_n})$.
	For any $\La\in\sI$, taking $n$ sufficiently large so that $\La\subset\La_n$,
	we see that $f^\La=\pi^{\La}f_n = \pi^{\La}(\pi^{\La_n} f_\infty)=\pi^{\La}f_\infty$
	by the tower property \eqref{eq: tower},
	which shows that $(f^\La)$ is the image of $f_\infty$ by \eqref{eq: nh} as desired.
	This gives our first assertion.
	
	Next, consider $f\in L^2(\mu)$ and the associated $(\pi^\La f)\in C_{L^2}(S^X)$.
	Again, Theorem \ref{thm: convergence} shows that
	the functions $f_n=\pi^{\La_n} f\in C(S^{\La_n})\subset C_\loc(S^X)$ for $n\in\bbN$
	strongly converge to $f_\infty\in L^2(\mu)$.
	By the previous argument, 
	since $\pi^\La f_\infty=\pi^\La f$ for any $\La\in\sI$, we see that the 
	image of $f$ and $f_\infty$ coincide in $C_{L^2}(S^X)$.
	Hence we have 
	\[
		\displaystyle\lim_{n\rightarrow\infty}f_n=f
	\] 
	in $C_{L^2}(S^X)$, which proves that $C_\loc(S^X)$ is dense in $C_{L^2}(S^X)$ as desired.
\end{proof}

\begin{rem}
	The space $C_{L^2}(S^X)$ is usually referred to as the \textit{space of Martingales bounded in $L^2$}
	(see for example \cite[Chapter 12]{Wil91}).
	The inclusions
	\[
		C_\loc(S^X) \subset C_{L^2}(S^X)\subset C_\col(S^X)
	\]
	give what is know as a \textit{Gelfand triple}.
\end{rem}

Next, we define the space of $L^2$-forms.  Assume that $\mu=\nu^{\otimes X}$ for some
probability measure $\nu$ supported on $S$.
For any $\La\in\sI$, we equip the space $\TS^\La$ with a norm given by
\[
	\dabs{\omega}^2_\mu\coloneqq\frac{1}{|E_\La|}\sum_{e\in E_\La}\dabs{\omega_e}^2_\mu
\]
for any $\omega=(\omega_e)\in\prod_{e\in E_\La}C(S^\La)$, where $\dabs{\omega_e}^2_\mu\coloneqq
E_\mu[\omega_e^2]$.     Following Theorem \ref{thm: convergence}, we define the space of $L^2$-forms in
\[
	C^1_\col(S^X)=\varprojlim_{\sI}C^1(S^\La)
\]
 as follows.

\begin{dfn}
	We define the space $C^1_{L^2}(S^X)$ of $L^2$-forms on $S^X$ by
	\[
		C^1_{L^2}(S^X)\coloneqq\bigl\{(\omega^\La)\in C^1_\col(S^X)\mid\sup_{\La\in\sI}\dabs{\omega^\La}_\mu<\infty\bigr\}\subset C^1_\col(S^X).
	\]
\end{dfn}

As in \S \ref{sec: Varadhan},
consider an action of a group $G$ on the locale $X$.  
We again assume that $G$ is torsion free
and that the action of $G$ on $X$ is free.  We let
\[
	\mathscr{C}\coloneqq C^1_{L^2}(S^X)\cap Z^1_\col(S^X)^G
	\subset  Z^1_\col(S^X)^G
\]
be the space of $G$-invariant closed $L^2$-forms, and we let $\mathscr{E}\coloneqq\overline{\partial(C_\unif(S^X)^G)}$ be the space of 
exact forms, where the bar denotes the closure of $\partial(C_\unif(S^X)^G)$ in $\mathscr{C}$
for the topology induced from the norm $\dabs{(\omega^\La)}_\mu\coloneqq
\sup_{\La\in\sI}\dabs{\omega^\La}_\mu$ for any $(\omega^\La)\in C^1_{L^2}(S^X)$.

By definition of the differential $\partial$,
for any $\xi\in\C{S}$, the function $\xi_X\in C^0_\unif(S^X)\subset C^0_\col(S^X)$
satisfies $\partial \xi_X=0$.  This shows that we have a natural inclusion 
$\C{S}\hookrightarrow H^0_\col(S^X)$ given by mapping $\xi$ to $\xi_X$.
Assume that $(X,E)$ is locale with a free action of a finitely generated abelian group $G$
of rank $d$,
such that the quotient $X/G$ is finite.  Such $(X,E)$ is a \emph{topological crystal}
in the sense of \cite[\S6.3]{Sun13}.  
We say that an interaction $\phi$ is \emph{simple}, if $c_\phi=1$ and the monoid
generated by $\xi(S)\subset\R$ is isomorphic to $\bbN$ or $\bbZ$.
We conjecture the following.

\begin{con}\label{con}
	Assume that the interaction $\phi$ is irreducibly quantified,
	and simple if $d=1$.  We assume in addition that the model satisfies a certain 
	spectral gap condition.
	Then the isomorphism of 
	Corollary \ref{crl: main} induces the isomorphism
	\[
		\mathscr{C}/\mathscr{E}\cong\Hom_\bbZ(G,\C{S}).
	\]
	In particular, we have a decomposition
	\[
		\mathscr{C}\cong\mathscr{E}\oplus\bigoplus_{j=1}^d\C{S}.
	\]
\end{con}

The above conjecture is the decomposition which is at the heart of the method proposed
by Varadhan to prove the hydrodynamic limit in the non-gradient case.  For the case that
the locale $(X,E)$ is the Euclidean lattice $(\bbZ^d,\bbE^d)$, the group $G=\bbZ^d$
with action on $(X,E)$ by translation, the set of states $S=\{0,1\}$, the interaction
is given by $\phi(s_1,s_2)=(s_2,s_1)$ for $(s_1,s_2)\in S\times S$, and $\nu$ is the probability 
measure on $S$
given by $\nu(\{1\})=p$ for $0<p<1$, Conjecture \ref{con} was proved by Funaki, Uchiyama and Yau \cite{FUY96}.  
We give in \cite{b} a proof of Conjecture \ref{con} for the Euclidean lattice $(\bbZ^d,\bbE^d)$
and group $G=\bbZ^d$.

\subsubsection*{Acknowledgment}

The authors would like to sincerely thank the organizers Iwao Kimura, Shinichi Kobayashi and Takashi Hara 
for the opportunity to give a talk at the RIMS Workshop Algebraic Number Theory and Related Topics 2020.
Giving a talk concerning probability theory at a workshop for algebraic number theory was a
very enriching experience.   We hope the application of algebraic and cohomological methods to  
stochastic models was interesting for number theorists, and look forward to further 
interactions and collaborations between various mathematical fields.
The authors are very grateful to Megumi Harada for carefully reading a preliminary version of this manuscript.
The authors would also like to thank the referee for very careful reading of the article and
detailed comments which helped to greatly improve the quality of the article.


\end{document}